\documentclass[11pt]{amsart}
\usepackage{amsthm}
\usepackage{amssymb}
\usepackage{latexsym}
\usepackage{amsfonts}
\usepackage{amsmath}
\usepackage{amstext}
\usepackage{amsrefs}
\newtheorem{theorem}{Theorem}[section]

\newtheorem{corollary}[theorem]{Corollary}

\newtheorem{definition}[theorem]{Definition}

\begin{document}

\title{A Note On The Isoperimetric Inequality And Its Stability}

\author{Xiang Gao}


\address{Department of Mathematics East China Normal University,~
Lane 500,~DongChuan Road,~Shanghai,~200241,~People's Republic of
China.}

\email{gaoxiangshuli@yahoo.cn}

\date{February 1,~2011}

\begin{abstract}
In this paper,~we deals with isoperimetric-type inequalities for
closed convex curves in the Euclidean plane~$ \mathbb{R}^2 $.~We
derive a family of parametric inequalities involving the following
geometric functionals associated to a given convex curve with a
simple Fourier series proof:~length,~area of the region included by
the curve,~area of the domain enclosed by the locus of curvature
centers and integral of the radius of curvature.~By using our
isoperimetric-type inequalities,~we also obtain some new geometric
Bonnesen-type inequalities.~Furthermore we investigate stability
properties of such inequalities~(near equality implies curve nearly
circular).
\end{abstract}

\keywords{isoperimetric inequality,~Fourier series,~stability}
\subjclass[2000]{~Primary 52A38;~Secondary 52A40}
\maketitle

\section{Introduction}

The classical isoperimetric inequality in the Euclidean plane~$
\mathbb{R}^2 $~states that:
\begin{theorem}[Isoperimetric Inequality]
If~$ \gamma$~is a simple closed curve of length L,~enclosing a
region of area A,~then
$$
L^2  - 4\pi A \ge 0, \eqno(1)
$$
and the equality holds if and only if~$ \gamma$~is a circle.
\end{theorem}
This fact was known to the ancient Greeks,~and the first
mathematical proof was only given in the 19th century by
Steiner.~Since then,~there have been many new proofs,~sharpened
forms, ~generalizations,~and applications of this famous inequality.

Suppose that~$ \gamma$~is a~$ \mathcal {C}_+^2$~closed and strictly
convex curve in the Euclidean plane~$ \mathbb{R}^2 $~with length
\emph{L},~area of the region included by the curve \emph{A},~and
area of the domain enclosed by the locus of curvature centers~$
\tilde A $.~Then there are also some interesting reverse
isoperimetric inequalities,~such as the inequality
$$
L^2  \le 4\pi \left( {A + | {\tilde A}|} \right), \eqno(2)
$$
proved by S.~L.~Pan and H.~Zhang in~[1],~and the
inequality~(3)~proved by S.~L.~Pan and J.~N.~Yang in~[2]:
$$
\int_0^{2\pi } {\rho \left( \theta  \right)} ^2 d\theta  \ge
\frac{{L^2  - 2\pi A}}{\pi }, \eqno(3)
$$
where~$ \rho$~is the radius of curvature and~$ \theta $~is the angle
between \emph{x}-axis and the outward normal vector at the
corresponding point \emph{p}.~Moreover the equalities
in~(2)~and~(3)~hold if and only if~$ \gamma$~is a circle.

It is obvious that if~$ \gamma$~is a circle,~then the locus of its
curvature centers is only a point,~and thus its area~$ \tilde A=0
$.~Conversely,~if~$ \tilde A=0 $,~then from the classical
isoperimetric inequality~(1)~and the reverse isoperimetric
inequality~(2),~it follows that the area \emph{A} and the length
\emph{L} of~$ \gamma$~satisfy~$ L^2  = 4\pi A$,~which implies that~$
\gamma$~is a circle,~and therefore the locus of curvature centers
of~$ \gamma$~is a point.

In this paper we deal with a family of parametric isoperimetric-type
inequalities for closed convex plane curves,~which is actually an
improved version of the reverse isoperimetric
inequalities~(2)~and~(3),~and one of the main results is as follows:
\begin{theorem}[Main Theorem]
Let~$ \gamma$~be a~$ \mathcal {C}_ + ^2$~closed and strictly convex
curve in the Euclidean plane~$ \mathbb{R}^2 $~with length L and
enclosing an area A,~then for arbitrary constants~$ \alpha , \beta ,
\lambda  , \delta$~satisfying
$$
\left\{ {\begin{array}{*{20}c}
   {2\alpha  + \delta  \ge 0}  \\
   {2\alpha  + 4\pi \beta  + \lambda  \ge 0}  \\
   {6\alpha  - \lambda  + 4\delta  \ge 0} , \\
\end{array}} \right. \eqno(4)
$$
we have
$$
\alpha \int_0^{2\pi } {\rho \left( \theta  \right)^2 } d\theta  +
\beta L^2  + \lambda A + \delta | {\tilde A} | \ge 0, \eqno(5)
$$
where~$ \rho$~is the curvature radius of~$ \gamma$~and~${\tilde A}
$~is the area of the domain enclosed by the locus of curvature
centers.~The equality holds if~$ \gamma$~is a circle and the
parameters~$ \alpha , \beta , \lambda  , \delta$~satisfy
$$
{2\alpha  + 4\pi \beta  + \lambda =0 }. \eqno(6)
$$
Moreover if the equality in~\emph{(5)}~holds and the parameters~$
\alpha , \beta , \lambda  , \delta$~satisfy
$$
\left\{ {\begin{array}{*{20}c}
   {2\alpha  + \delta  > 0}  \\
   {2\alpha  + 4\pi \beta  + \lambda  = 0}  \\
   {6\alpha  - \lambda  + 4\delta  = 0},  \\
\end{array}} \right.\eqno(7)
$$
then the Minkowski support function of~$ \gamma$~is of the form~$
p\left( \theta \right) = a_0 + a_1 \cos \theta  + b_1 \sin \theta  +
a_2 \cos 2\theta  + b_2 \sin 2\theta$.
\end{theorem}
{\it Remark 1.}~When~$ \alpha  = 0,\beta  =  - 1,\lambda  = \delta =
4\pi$,~(4)~satisfies and the isoperimetric inequality~(5)~turns
into~(2).~When~$ \alpha  = 1,\beta  =  - \frac{1}{\pi },\lambda  =
2,\delta  = 0$,~(4)~also satisfies and we
obtain~(3).~Hence~(5)~could also be regarded as a reverse
isoperimetric inequality.~Furthermore,~if we select other values of
the parameters~$ \alpha , \beta , \lambda  ,
\delta$~satisfying~(4),~we can obtain some new geometric
Bonnesen-type inequalities~[3]:
\begin{corollary}
Let~$ \gamma$~be a~$ \mathcal {C}_ + ^2$~closed and strictly convex
curve in the Euclidean plane~$ \mathbb{R}^2 $~with length L and
enclosing an area A,~we have
$$
L^2  \le 4\pi A + \pi |\tilde A|,\eqno(8)
$$
$$
\int_0^{2\pi } {\rho \left( \theta  \right)^2 } d\theta  \ge
\frac{{L^2 }}{\pi } - 2A + |\tilde A|,\eqno(9)
$$
and
$$
\mathop {\max }\limits_{\theta  \in \left[ {0,2\pi } \right]} \rho
\left( \theta  \right)^2  \ge \frac{{1 }}{2 \pi }  \left( \frac{{L^2
}}{\pi } - 2A + |\tilde A| \right), \eqno(10)
$$
where~$ \rho$~is the curvature radius of~$ \gamma$~and~${\tilde A}
$~is the area of the domain enclosed by the locus of curvature
centers.~Furthermore,~\emph{(8)}~and~\emph{(9)}~are improved
versions of~\emph{(2)}~and~\emph{(3)},~and the equalities
in~\emph{(8)},~\emph{(9)}~and~\emph{(10)}~hold if~$ \gamma$~is a
circle.~Moreover if the equalities
in~\emph{(8)},~\emph{(9)}~and~\emph{(10)}~hold,~then the Minkowski
support function of~$ \gamma$~is of the form~$ p\left( \theta
\right) = a_0 + a_1 \cos \theta  + b_1 \sin \theta  + a_2 \cos
2\theta  + b_2 \sin 2\theta$.
\end{corollary}

The stability problem associated with isoperimetric inequality is
also interesting and significant.

Recently in [4], S. L. Pan and H. P. Xu obtained the following
stability estimates for the reverse isoperimetric inequality~(2)~by
comparing a convex body \emph{K} with its Steiner disk.
\[
\begin{split}
 h_1 \left( {K,S\left( K \right)} \right)^2  &= \left( {\mathop {\max }\limits_u \left| {p_K \left( u \right) - p_{S\left( K \right)} \left( u \right)} \right|} \right)^2  \\
  & \le \frac{{4\pi ^2  - 33}}{{96\pi ^2 }}\left( {4\pi \left( {A\left( K \right) + | {\tilde A\left( K \right)} |} \right) - L^2 \left( K \right)} \right), \\
 \end{split}
\]
\[
\begin{split}
 h_2 \left( {K,S\left( K \right)} \right)^2 & =  \int_0^{2\pi } {\left| {p_K \left( \theta  \right) - p_{S\left( K \right)} \left( \theta  \right)} \right|} ^2 d\theta  \\
  & \le  \frac{1}{{18\pi }}\left( {4\pi \left( {A\left( K \right) + | {\tilde A\left( K \right)} |} \right) - L^2 \left( K \right)} \right), \\
 \end{split}
\]
where~$ p_K \left( \theta  \right)$~denotes the Minkowski support
function of a given convex body \emph{K},~and~$S\left( K
\right)$~denotes the Steiner disc associated with \emph{K}~(see
section 4 for the definition)~which satisfies
$$
4\pi \left( {A\left( {S\left( K \right)} \right) + | {\tilde A\left(
{S\left( K \right)} \right)} |} \right) - L^2 \left( {S\left( K
\right)} \right) = 0.\eqno(11)
$$

For arbitrary~$ \varepsilon >0 $~such that~$ \varphi \left( K
\right) = 4\pi \left( {A\left( K \right) + | {\tilde A\left( K
\right)} |} \right) - L^2 \left( K \right) < \varepsilon$,~by the
stability estimates for inequality above and~(11)~it follows that
$$
\max \left\{ {h_1 \left( {K,S\left( K \right)} \right)^2,h_2 \left(
{K,S\left( K \right)} \right)^2} \right\} \le C| {\varphi \left( K
\right) - \varphi \left( {S\left( K \right)} \right)} | <
C\varepsilon,
$$
which implies that the reverse isoperimetric inequality~(2)~does
have a good stability behaviour with respect to both Hausdorff
distance and~$ L^2 $-metric.

The paper is organized as follows.~In section 2,~we recall some
basic facts about plane convex geometry.~In section 3,~we provide a
simpler proof of Theorem 1.2 by using Fourier series,~which is
different from the approach in~[1]~and~[2].~In section 4,~we
investigate stability properties of inequality~(5)~(near equality
implies curve nearly circular).~We believe that our trick could be
used to derive more interesting isoperimetric inequalities.

\section{Geometric Quantities and Their Fourier Series}

In this section,~we recall some basic facts about convex plane curve
which will be used later.~In this paper we always assume
that~$\gamma $~is a closed and convex plane curve which is
sufficiently regular,~actually it should be a~$ \mathcal {C}_ +
^2$~closed and strictly convex curve in the plane~$ \mathbb{R}^2
$,~such that the radius of curvature can be defined and the Fourier
series needed in the proof convergent uniformly.~The details can be
found in the classical literature~[5].

Let~$ p\left( \theta  \right)$~denote the Minkowski support function
of curve~$\gamma \left( \theta \right)$,~where~$ \theta $~is the
angle between \emph{x}-axis and the outward normal vector at the
corresponding point \emph{p}.~It gives us the parametrization
of~$\gamma \left( \theta  \right)$~in terms of~$ \theta$~as follows:
$$
\gamma \left( \theta  \right) = \left( {\gamma _1 \left( \theta
\right),\gamma _2 \left( \theta  \right)} \right) = \left( {p\left(
\theta  \right)\cos \theta  - p'\left( \theta  \right)\sin \theta
,p\left( \theta  \right)\sin \theta  + p'\left( \theta  \right)\cos
\theta } \right).
$$
Therefore the curvature~$ k\left( \theta  \right) $~and the radius
of curvature~$ \rho\left( \theta  \right) $~of~$\gamma \left( \theta
\right)$~can be calculated by
$$
k\left( \theta  \right) = \frac{{d\theta }}{{ds}} =
\frac{1}{{p\left( \theta  \right) + p''\left( \theta  \right)}} >
0
$$
and
$$
\rho \left( \theta  \right) = \frac{{ds}}{{d\theta }} = p\left(
\theta  \right) + p''\left( \theta  \right) > 0.
$$
The length \emph{L} of~$\gamma \left( \theta  \right)$~and the area
\emph{A} it bounds can be also calculated respectively by
$$
L = \int_\gamma  {ds = } \int_0^{2\pi } {p \left( \theta \right)}
d\theta
$$
and
$$
A = \frac{1}{2}\int_\gamma  {p\left( \theta \right)ds = }
\frac{1}{2}\int_0^{2\pi } {\left( {p\left( \theta \right)^2  -
p'\left( \theta  \right)^2 } \right)} d\theta.
$$
At the same time,~we could obtain the locus of centers of curvature
of~$\gamma \left( \theta  \right)$~as follow
$$
\beta \left( \theta  \right) = \gamma \left( \theta  \right) + \rho
\left( \theta  \right)N\left( \theta  \right) = \left( { - p'\left(
\theta  \right)\sin \theta  - p''\left( \theta  \right)\cos \theta
,p'\left( \theta  \right)\cos \theta  - p''\left( \theta \right)\sin
\theta } \right),
$$
and the oriented area of the domain enclosed by~$ \beta \left(
\theta \right)$~is given by
$$
\tilde A = \frac{1}{2}\int_0^{2\pi } {\left( {p'\left( \theta
\right)^2  - p''\left( \theta  \right)^2 } \right)} d\theta.
$$

Since the Minkowski support function of a given convex body \emph{K}
is always continuous,~bounded and~$ 2\pi$-periodic,~it has a Fourier
series of the form
$$
p\left( \theta  \right) = a_0  + \sum\limits_{n = 1}^\infty  {\left(
{a_n \cos n\theta  + b_n \sin n\theta } \right)}. \eqno(12)
$$
Differentiation of~(12)~with respect to~$ \theta$~gives us
$$
p'\left( \theta  \right) = \sum\limits_{n = 1}^\infty  {n\left( { -
a_n \sin n\theta  + b_n \cos n\theta } \right)} \eqno(13)
$$
and
$$
p''\left( \theta \right) =  - \sum\limits_{n = 1}^\infty  {n^2
\left( {a_n \cos n\theta  + b_n \sin n\theta } \right)}. \eqno(14)
$$
Thus by~(12),~(13),~(14)~and the Parseval equality we could express
these geometric quantities in terms of the Fourier coefficients of~$
p\left( \theta  \right)$
$$
\begin{array}{l}
 \rho \left( \theta  \right) = p\left( \theta  \right) + p''\left( \theta  \right) \\
  \qquad = a_0  + \sum\limits_{n = 1}^\infty  {\left( {a_n \cos n\theta  + b_n \sin n\theta } \right)}  - \sum\limits_{n = 1}^\infty  {n^2 \left( {a_n \cos n\theta  + b_n \sin n\theta } \right)},  \\
 \end{array} \eqno(15)
$$
$$
L\left( K \right) = 2\pi a_0, \eqno(16)
$$
$$
A = \pi a_0^2  - \frac{\pi }{2}\sum\limits_{n = 2}^\infty  {\left(
{n^2  - 1} \right)\left( {a_n^2  + b_n^2 } \right)} ,\eqno(17)
$$
$$
| {\tilde A} | = \frac{\pi }{2}\sum\limits_{n = 2}^\infty  {n^2
\left( {n^2  - 1} \right)\left( {a_n^2  + b_n^2 } \right)}.
\eqno(18)
$$

\section{Proof of The Main Theorems}

\begin{proof}[Proof of Theorem 1.2]
Firstly from~(15),~one can easily get
$$
\begin{array}{l}
 \quad \int_0^{2\pi } {\rho \left( \theta  \right)^2 } d\theta  \\
\\
  = 2\left( {\pi a_0^2  - \frac{\pi }{2}\sum\limits_{n = 2}^\infty  {\left( {n^2  - 1} \right)\left( {a_n^2  + b_n^2 } \right)}  + \frac{\pi }{2}\sum\limits_{n = 2}^\infty  {n^2 \left( {n^2  - 1} \right)\left( {a_n^2  + b_n^2 } \right)} } \right) \\
  = 2\pi \left( {a_0^2  + \frac{1}{2}\sum\limits_{n = 2}^\infty  {\left( {n^2  - 1} \right)^2 \left( {a_n^2  + b_n^2 } \right)} } \right), \\
 \end{array}
$$
thus by using~(16),~(17)~and~(18)~we have
$$
\begin{array}{l}
 \quad \alpha \int_0^{2\pi } {\rho \left( \theta  \right)^2 } d\theta  + \beta L^2  + \lambda A + \delta | {\tilde A} | \\
\\
  = 2\pi \alpha \left( {a_0^2  + \frac{1}{2}\sum\limits_{n = 2}^\infty  {\left( {n^2  - 1} \right)^2 \left( {a_n^2  + b_n^2 } \right)} } \right) + \beta \left( {2\pi a_0 } \right)^2  \\
 \quad + \lambda \left( {\pi a_0^2  - \frac{\pi }{2}\sum\limits_{n = 2}^\infty  {\left( {n^2  - 1} \right)\left( {a_n^2  + b_n^2 } \right)} } \right) + \frac{\pi }{2}\delta \sum\limits_{n = 2}^\infty  {n^2 \left( {n^2  - 1} \right)\left( {a_n^2  + b_n^2 } \right)}  \\
  = \pi a_0^2 \left( {2\alpha  + 4\pi \beta  + \lambda } \right) + \frac{\pi }{2}\sum\limits_{n = 2}^\infty  {\left( {2\alpha \left( {n^2  - 1} \right) - \lambda  + \delta n^2 } \right)\left( {n^2  - 1} \right)\left( {a_n^2  + b_n^2 } \right)}  \\
 \end{array} \eqno(19)
$$
It follows from~(4)~that~(19)~is nonnegative,~which completes the
proof of inequality~(5).

Furthermore,~if~$ \gamma$~is a circle,~by the equality conditions
in~(2)~and~(3)~we have
$$
L^2  = 4\pi \left( {A + |\tilde A|} \right) = 4\pi A
$$
and
$$
\int_0^{2\pi } {\rho \left( \theta  \right)^2 } d\theta  =
\frac{{L^2  - 2\pi A}}{\pi } = 2A.
$$
Hence
$$
\begin{array}{l}
 \quad \alpha \int_0^{2\pi } {\rho \left( \theta  \right)^2 } d\theta  + \beta L^2  + \lambda A + \delta |\tilde A| \\
\\
  = 2\alpha A + 4\pi \beta A + \lambda A \\
\\
  = \left( {2\alpha  + 4\pi \beta  + \lambda } \right)A \\
 \end{array}
$$
then for the parameters~$ \alpha , \beta , \lambda  ,
\delta$~satisfying~(6)~we have
$$
\alpha \int_0^{2\pi } {\rho \left( \theta  \right)^2 } d\theta  +
\beta L^2  + \lambda A + \delta |\tilde A| = 0
$$
On the other hand,~if equality holds in~(5):
\[
\begin{split}
 0 &= \alpha \int_0^{2\pi } {\rho \left( \theta  \right)^2 } d\theta  + \beta L^2  + \lambda A + \delta |\tilde A| \\
   &= \pi a_0^2 \left( {2\alpha  + 4\pi \beta  + \lambda } \right) + \frac{{3\pi }}{2}\left( {6\alpha  - \lambda  + 4\delta } \right)\left( {a_2^2  + b_2^2 } \right) \\
   & \quad + \frac{\pi }{2}\sum\limits_{n = 3}^\infty  {\left( {2\alpha \left( {n^2  - 1} \right) - \lambda  + \delta n^2 } \right)\left( {n^2  - 1} \right)\left( {a_n^2  + b_n^2 } \right)}  \\
 \end{split}
\]
then by the condition~(7):
$$
\left\{ {\begin{array}{*{20}c}
   {2\alpha  + \delta  > 0}  \\
   {2\alpha  + 4\pi \beta  + \lambda  = 0}  \\
   {6\alpha  - \lambda  + 4\delta  = 0}  \\
\end{array}} \right.
$$
we have
$$
0 = \frac{\pi }{2}\sum\limits_{n = 3}^\infty  {\left( {2\alpha
\left( {n^2  - 1} \right) - \lambda  + \delta n^2 } \right)\left(
{n^2  - 1} \right)\left( {a_n^2  + b_n^2 } \right)}
$$
and
$$
\left( {2\alpha \left( {n^2  - 1} \right) - \lambda  + \delta n^2 }
\right)\left( {n^2  - 1} \right) > 0
$$
for~$ n \ge 3$.~Thus~$ a_n  = b_n  = 0$~for~$ n \ge 3$~and the
Minkowski support function of~$ \gamma$~is of the form~$ p\left(
\theta \right) = a_0 + a_1 \cos \theta  + b_1 \sin \theta  + a_2
\cos 2\theta  + b_2 \sin 2\theta$. This completes the proof of
Theorem 1.2.
\end{proof}
\begin{proof}[Proof of Corollary 1.3]
Let~$ \alpha  = 0,\beta  =  - 1,\lambda  =  4\pi, \delta  = \pi$~we
obtain~(8),~and let~$ \alpha = 1,\beta  =  - \frac{1}{\pi },\lambda
= 2,\delta  = - 1$,~we can derive~(9).~Moreover the equality
conditions in~(8)~and~(9)~follows directly from the equality
conditions in~(5).

On the other hand,~inequality~(10)~is an easy consequence
of~(9),~and if~$ \gamma$~is a circle,~the equality holds
directly.~Conversely,~since
$$
\mathop {\max }\limits_{\theta  \in \left[ {0,2\pi } \right]} \rho
\left( \theta  \right)^2  \ge \frac{1}{2\pi }\int_0^{2\pi } {\rho
\left( \theta  \right)^2 } d\theta  \ge \frac{1}{2\pi }\left(
{\frac{{L^2 }}{\pi } - 2A + |\tilde A|} \right),
$$
if equality holds in~(10),~we have
$$
\int_0^{2\pi } {\rho \left( \theta  \right)^2 } d\theta  =
\frac{{L^2 }}{\pi } - 2A + |\tilde A|.
$$
By the equality condition of~(9),~it follows that the Minkowski
support function of~$ \gamma$~is of the form~$ p\left( \theta
\right) = a_0 + a_1 \cos \theta  + b_1 \sin \theta + a_2 \cos
2\theta  + b_2 \sin 2\theta$.
\end{proof}

\section{The Stability Property of The Isoperimetric Inequality}

Let \emph{K} and \emph{M} be two convex bodies with respective
Minkowski support functions~$ p_K$~and~$ p_M$.~The most frequently
used function to measure the deviation between \emph{K} and \emph{M}
is the Hausdorff distance:
$$
h_1 \left( {K,M} \right) = \mathop {\max }\limits_u | {p_K \left( u
\right) - p_M \left( u \right)} |.
$$
Another distance is defined by means of the~$ L^2 $-norm of the
support functions,~that is
$$
h_2 \left( {K,M} \right) = \left( {\int_0^{2\pi } {| {p_K \left(
\theta  \right) - p_M \left( \theta  \right)} |} ^2 d\theta }
\right)^{\frac{1}{2}},
$$
where~$ \theta $~is the angle between \emph{x}-axis and the outward
normal vector at the corresponding point \emph{p}.~It is obvious
that~$ h_1 \left( {K,M} \right) = 0$~or~$ h_2 \left( {K,M} \right) =
0$~if and only if \emph{K = M}.

We also recall the definition of Steiner disc \emph{S(K)} of a
planar convex body \emph{K}.
\begin{definition}
The Steiner disc of a convex body K,~denoted by S(K) is the circular
disc with radius~$\frac{{L\left( K \right)}}{{2\pi }}$~and center at
the Steiner point~$ \overrightarrow s \left( K \right) $~which can
be defined in terms of the Minkowski support function~$ p_K \left(
\theta \right)$:
$$
\overrightarrow s \left( K \right) = \frac{1}{\pi }\int_0^{2\pi } {
\overrightarrow u \left( \theta  \right)p_K \left( \theta  \right)}
d\theta,
$$
where~$ \overrightarrow u \left( \theta  \right) $~is a unit tangent
vector at the corresponding point p,~and~$ L\left( K \right)
$~denotes the perimeter of the domain K.
\end{definition}

We now derive a stability version of~(5)~with respect to both
Hausdorff distance~$ h_1 $~and~$ h_2 $~metric.
\begin{theorem}
Let K be a domain enclosed by a~$ \mathcal {C}_+^2$~closed and
strictly convex plane curve~$ \gamma$~with area A(K) and perimeter
L(K),~and let~$ \tilde A (K) $~denote the oriented area of the
domain enclosed by the locus of curvature centers of~$ \gamma$,~S(K)
denotes the Steiner disc associated with K.~Then for arbitrary
constants~$ \alpha , \beta , \lambda  , \delta$~which satisfy
$$
\left\{ {\begin{array}{*{20}c}
   {2\alpha  + \delta  \ge 0}  \\
   {2\alpha  + 4\pi \beta  + \lambda  \ge 0}  \\
   {6\alpha  - \lambda  + 4\delta  > 0} , \\
\end{array}} \right.  \eqno(20)
$$
we have
$$
h_1 \left( {K,S\left( K \right)} \right)^2  \le C \left( \alpha ,
 \lambda  , \delta \right) \left( {\alpha \int_0^{2\pi }
{\rho \left( \theta  \right)^2 } d\theta  + \beta L^2 \left( K
\right) + \lambda A \left( K \right) + \delta | {\tilde A\left( K
\right)} |} \right), \eqno(21)
$$
where~$ C \left( \alpha ,
 \lambda  , \delta \right) = \max \left\{ {1,\frac{2}{\pi}\sum\limits_{n = 2}^\infty
{\frac{1}{{\left( {2\alpha \left( {n^2  - 1} \right) - \lambda  +
\delta n^2 } \right)\left( {n^2  - 1} \right)}}} } \right\} $.~The
equality holds if~$ \gamma$~is a circle and the parameters~$ \alpha
, \beta , \lambda  , \delta$~satisfy
$$
{2\alpha  + 4\pi \beta  + \lambda =0 }.
$$
\end{theorem}
\begin{proof}
We may assume~$ \overrightarrow s \left( K \right) = 0$,~because
of~(12)~and~(16),~the support functions~$ p_K$~and~$ p_{S(K)}$~have
the following Fourier series:
$$
p_K \left( \theta  \right) = \frac{{L\left( K \right)}}{{2\pi }} +
\sum\limits_{n = 2}^\infty  {\left( {a_n \cos n\theta  + b_n \sin
n\theta } \right)}\eqno(22)
$$
and
$$
p_{S\left( K \right)} \left( \theta  \right) = \frac{{L\left( K
\right)}}{{2\pi }}. \eqno(23)
$$
One can observe that~(22)~and~(23)~yield an explicit expression~(in
terms of the Fourier coefficients)~for the quantity:
$$
\begin{array}{l}
 \quad \alpha \int_0^{2\pi } {\rho \left( \theta  \right)^2 } d\theta  + \beta L^2 \left( K \right) + \lambda A \left( K \right) + \delta | {\tilde A \left( K \right) } | \\
\\
  = \pi a_0^2 \left( {2\alpha  + 4\pi \beta  + \lambda } \right) + \frac{\pi }{2}\sum\limits_{n = 2}^\infty  {\left( {2\alpha \left( {n^2  - 1} \right) - \lambda  + \delta n^2 } \right)\left( {n^2  - 1} \right)\left( {a_n^2  + b_n^2 } \right)}.  \\
 \end{array} \eqno(24)
$$
Since it is easily seen that
$$
 | {a_n \cos n\theta  + b_n \sin
n\theta } | \le \sqrt {a_n^2  + b_n^2 },
$$
it follows that
\[
\begin{split}
 | {p_K \left( \theta  \right) - p_{S\left( K \right)} \left( \theta  \right)} | &= \left | {\frac{{L\left( K \right)}}{{2\pi }} + \sum\limits_{n = 2}^\infty  {\left( {a_n \cos n\theta  + b_n \sin n\theta } \right)}  - \frac{{L\left( K \right)}}{{2\pi }}} \right | \\
 & \le \sum\limits_{n = 2}^\infty  {| {a_n \cos n\theta  + b_n \sin n\theta } |} \\
 & \le \sum\limits_{n = 2}^\infty  {\sqrt {a_n^2  + b_n^2 } }.  \\
 \end{split}
\]
Using Holder's inequality,~together with~(24)~we have
\[
\begin{split}
 h_1 \left( {K,S\left( K \right)} \right)^2 & \le \left( {\sum\limits_{n = 2}^\infty  {\sqrt {a_n^2  + b_n^2 } } } \right)^2  \\
 & \le \pi a_0^2 \left( {2\alpha  + 4\pi \beta  + \lambda } \right) + \left( \frac{2}{\pi} {\sum\limits_{n = 2}^\infty  {\frac{1}{{\left( {2\alpha \left( {n^2  - 1} \right) - \lambda  + \delta n^2 } \right)\left( {n^2  - 1} \right)}}} }
 \right)\\
 & \quad \left( \frac{\pi }{2} {\sum\limits_{n = 2}^\infty  {\left( {2\alpha \left( {n^2  - 1} \right) - \lambda  + \delta n^2 } \right)\left( {n^2  - 1} \right)\left( {a_n^2  + b_n^2 } \right)} } \right) \\
 & \le \max \left\{ {1,\frac{2}{\pi} \sum\limits_{n = 2}^\infty  {\frac{1}{{\left(
{2\alpha \left( {n^2  - 1} \right) - \lambda  + \delta n^2 }
\right)\left( {n^2  - 1} \right)}}} } \right\}\\
 & \quad \left( {\alpha \int_0^{2\pi } {\rho \left( \theta  \right)^2 } d\theta  + \beta L^2 \left( K \right) + \lambda A \left( K \right) + \delta | {\tilde A \left( K \right)} |} \right), \\
\end{split}
\]
for arbitrary constants~$ \alpha , \beta , \lambda ,
\delta$~satisfying~(20).

Furthermore,~if~$ \gamma$~is a circle,~as the proof of Theorem 1.2
we have
$$
\begin{array}{l}
 \quad \alpha \int_0^{2\pi } {\rho \left( \theta  \right)^2 } d\theta  + \beta L^2 \left( K \right) + \lambda A \left( K \right) + \delta |\tilde A \left( K \right) | \\
\\
  = \left( {2\alpha  + 4\pi \beta  + \lambda } \right)A \\
 \end{array}
$$
If the parameters~$ \alpha , \beta , \lambda  , \delta$~satisfy~$
{2\alpha  + 4\pi \beta  + \lambda =0 } $,~then we have
$$
\alpha \int_0^{2\pi } {\rho \left( \theta  \right)^2 } d\theta  +
\beta L^2 \left( K \right) + \lambda A \left( K \right) + \delta
|\tilde A \left( K \right)| = 0
$$
It is obvious that~$ h_1 \left( {K,S(K)} \right) = 0$,~thus equality
holds in~(21).
\end{proof}
\begin{theorem}
Under the same assumptions of Theorem 4.2,~then for arbitrary
constants~$ \alpha , \beta , \lambda  , \delta$~which satisfy
$$
\left\{ {\begin{array}{*{20}c}
   {2\alpha  + \delta  \ge 0}  \\
   {2\alpha  + 4\pi \beta  + \lambda  \ge 0}  \\
   {18\alpha  - 3\lambda  + 12\delta  - 2 \ge 0},  \\
\end{array}} \right. \eqno(25)
$$
we have
$$
h_2 \left( {K,S\left( K \right)} \right)^2  \le \alpha \int_0^{2\pi
} {\rho \left( \theta  \right)^2 } d\theta  + \beta L^2  + \lambda A
+ \delta | {\tilde A} |.  \eqno(26)
$$
The equality holds if~$ \gamma$~is a circle and the parameters~$
\alpha , \beta , \lambda  , \delta$~satisfy
$$
{2\alpha  + 4\pi \beta  + \lambda =0 }.
$$
Moreover if the equality in~\emph{(26)}~holds and the parameters~$
\alpha , \beta , \lambda  , \delta$~satisfy
$$
\left\{ {\begin{array}{*{20}c}
   {2\alpha  + \delta  > 0}  \\
   {2\alpha  + 4\pi \beta  + \lambda  = 0}  \\
   {18\alpha  - 3\lambda  + 12\delta  - 2 = 0},   \\
\end{array}} \right.  \eqno(27)
$$
then the Minkowski support function of~$ \gamma$~is of the form~$
p\left( \theta \right) = a_0 + a_1 \cos \theta  + b_1 \sin \theta  +
a_2 \cos 2\theta  + b_2 \sin 2\theta$.
\end{theorem}
\begin{proof}
As the proof of Theorem 4.2,~we use Parseval's
equality,~(22)~and~(23)~to deduce that
$$
h_2 \left( {K,S\left( K \right)} \right)^2  = \int_0^{2\pi } {| {p_K
\left( \theta  \right) - p_{S\left( K \right)} \left( \theta
\right)} |} ^2 d\theta  = \pi \sum\limits_{n = 2}^\infty {\left(
{a_n^2  + b_n^2 } \right)},
$$
together with~(24)~one gets that
$$
\begin{array}{l}
 \quad \left( {\alpha \int_0^{2\pi } {\rho \left( \theta  \right)^2 } d\theta  + \beta L^2 \left( K \right) + \lambda A \left( K \right) + \delta | {\tilde A \left( K \right) } |} \right) - h_2 \left( {K,S\left( K \right)} \right)^2  \\
\\
  = \pi a_0^2 \left( {2\alpha  + 4\pi \beta  + \lambda } \right) + \frac{\pi }{2}\sum\limits_{n = 2}^\infty  {\left( {2\alpha \left( {n^2  - 1} \right) - \lambda  + \delta n^2 } \right)\left( {n^2  - 1} \right)\left( {a_n^2  + b_n^2 } \right)} \\
 \quad - \pi \sum\limits_{n = 2}^\infty  {\left( {a_n^2  + b_n^2 } \right)}  \\
  = \pi a_0^2 \left( {2\alpha  + 4\pi \beta  + \lambda } \right)
  + \frac{\pi }{2}\sum\limits_{n = 2}^\infty  {\left( {\left( {2\alpha \left( {n^2  - 1} \right) - \lambda  + \delta n^2 } \right)\left( {n^2  - 1} \right) - 2} \right)\left( {a_n^2  + b_n^2 } \right)}. \\
 \end{array}
$$
Hence for arbitrary constants~$ \alpha , \beta , \lambda ,
\delta$~satisfying~(25),~we have
$$
\quad \left( {\alpha \int_0^{2\pi } {\rho \left( \theta  \right)^2 }
d\theta  + \beta L^2 \left( K \right) + \lambda A \left( K \right) +
\delta | {\tilde A \left( K \right) } |} \right) - h_2 \left(
{K,S\left( K \right)} \right)^2  \ge 0,
$$
which implies the following stability result:
$$
h_2 \left( {K,S\left( K \right)} \right)^2  \le \alpha \int_0^{2\pi
} {\rho \left( \theta  \right)^2 } d\theta  + \beta L^2 \left( K
\right) + \lambda A \left( K \right)  + \delta | {\tilde A \left( K
\right) } |.
$$

Furthermore,~if~$ \gamma$~is a circle,~as the proof of Theorem
4.2,~we have equality in~(26).~Conversely,~if equality holds
in~(26):
\[
\begin{split}
 0 &= \left( {\alpha \int_0^{2\pi } {\rho \left( \theta  \right)^2 } d\theta  + \beta L^2 \left( K \right) + \lambda A\left( K \right) + \delta |\tilde A\left( K \right)|} \right) - h_2 \left( {K,S\left( K \right)} \right)^2  \\
   &= \pi a_0^2 \left( {2\alpha  + 4\pi \beta  + \lambda } \right) + \frac{\pi }{2}\sum\limits_{n = 2}^\infty  {\left( {\left( {2\alpha \left( {n^2  - 1} \right) - \lambda  + \delta n^2 } \right)\left( {n^2  - 1} \right) - 2} \right)\left( {a_n^2  + b_n^2 } \right)}  \\
 \end{split}
\]
then by the condition~(27),~we have
$$
0 = \frac{\pi }{2}\sum\limits_{n = 3}^\infty  {\left( {\left(
{2\alpha \left( {n^2  - 1} \right) - \lambda  + \delta n^2 }
\right)\left( {n^2  - 1} \right) - 2} \right)\left( {a_n^2  + b_n^2
} \right)}
$$
and
$$
\left( {2\alpha \left( {n^2  - 1} \right) - \lambda  + \delta n^2 }
\right)\left( {n^2  - 1} \right) - 2 > 0
$$
for~$ n \ge 3$.~Thus~$ a_n  = b_n  = 0$~for~$ n \ge 3$~and the
Minkowski support function of~$ \gamma$~is of the form~$ p\left(
\theta \right) = a_0 + a_1 \cos \theta  + b_1 \sin \theta  + a_2
\cos 2\theta  + b_2 \sin 2\theta$. This completes the proof of
Theorem 4.3.
\end{proof}
{\it Remark 2.}~The combination of Theorem 4.2 and 4.3 leads to
$$
\begin{array}{l}
 \quad \max \left\{ {h_1 \left( {K,S\left( K \right)} \right)^2 ,h_2 \left( {K,S\left( K \right)} \right)^2 } \right\} \\
\\
  \le C\left( \alpha ,
 \lambda  , \delta \right)\left( {\alpha \int_0^{2\pi } {\rho \left( \theta  \right)^2 } d\theta  + \beta L^2 \left( K \right) + \lambda A \left( K \right)  + \delta | {\tilde A \left( K \right) } |} \right), \\
 \end{array} \eqno(30)
$$
where~$ C \left( \alpha ,
 \lambda  , \delta \right) = \max \left\{ {1,\frac{2}{\pi}\sum\limits_{n = 2}^\infty
{\frac{1}{{\left( {2\alpha \left( {n^2  - 1} \right) - \lambda  +
\delta n^2 } \right)\left( {n^2  - 1} \right)}}} } \right\} $,~which
states that the isoperimetric inequality~(5)~does have a good
stability behaviour with respect to both Hausdorff distance and~$
L^2 $-metric.

\section{Acknowledgment}
I would especially like to express my appreciation to my advisor
professor Yu Zheng for longtime encouragement and meaningful
discussions.~I would also especially like to thank the referee for
meaningful suggestions that led to improvement of the article.


\begin{thebibliography}{99}

\bibitem{1}
S.~L.~Pan,~H.~Zhang,~\emph{A reverse isoperimetric inequality for
convex plane curves},\\Beitrage Algebra Geom.~48~(2007)~303-308.

\bibitem{2}
S.~L.~Pan,~J.~N.~Yang,~\emph{On a non-local perimeter-preserving
curve evolution problem for convex plane curves},~Manuscripta
Math.~127,~(2008)~469-484.

\bibitem{3}
R.~Schneider,~\emph{Convex bodies: The Brunn-Minkowski theory},
~Encyclopedia of Mathematics and its Applications,~44.~Cambridge
University Press,~Cambridge,~1993.

\bibitem{4}
S.~L.~Pan,~H.~P.~Xu,~\emph{Stability of a reverse isoperimetric
inequality},~J.~Math.~Anal.\\~Appl.~350~(2009)~348-353.

\bibitem{5}
H.~Groemer,~\emph{Geometric applications of Fourier series and
spherical harmonics}, ~Encyclopedia of Mathematics and its
Applications,~61.~Cambridge University Press,~Cambridge,~1996.

\end{thebibliography}
\end{document}